\documentclass[12pt]{article}
\usepackage[utf8]{inputenc}
\usepackage{amsmath, amsfonts, amsthm,fullpage,doi}
\let\epsilon=\varepsilon
\title{Alternating paths in oriented graphs with large semidegree}
\author{Jozef Skokan \thanks{Department of Mathematics, London School of Economics, Houghton Street, London, WC2A 2AE.} \and Mykhaylo Tyomkyn \thanks{Department of Applied Mathematics, Charles University. Email: tyomkyn$@$kam.mff.cuni.cz.
Supported in part by ERC Synergy Grant DYNASNET 810115 and GA\v{C}R Grant 22-19073S.}}
\date{\today}

\newtheorem{theorem}{Theorem}[section]
\newtheorem{lemma}[theorem]{Lemma}
\newtheorem{prop}[theorem]{Proposition}
\newtheorem{cor}[theorem]{Corollary}
\newtheorem{conj}[theorem]{Conjecture}

\begin{document}

\maketitle

\begin{abstract}
In new progress on conjectures of Stein, 
and Addario-Berry, Havet, Linhares Sales, Reed and Thomass\'e, 
we prove that every oriented graph with all in- and out-degrees greater than $5k/8$ contains an alternating path of length $k$. This improves on previous results of Klimo\v{s}ov\'a and Stein, 
and Chen, Hou and Zhou. 
\end{abstract}

\section{Introduction}

A classical topic in graph theory is establishing minimum degree conditions for existence of long paths and cycles in graphs. This research direction has a rich history, going back to works of Dirac~\cite{Di}, P\'osa~\cite{Po} and others. The present paper is a part of the ongoing effort to establish analogous results in oriented graphs. 

An oriented graph is a tuple $G=(V,\vec{D})$ where $V$ is a finite set of \emph{vertices} and $\vec{D}$ is a set of \emph{(directed) edges}. Each edge is an ordered pair of distinct vertices. Moreover, for all $u$ and $v$ at most one of $(u,v)$ and $(v,u)$ is an edge. 
The \emph{in-degree} and \emph{out-degree} of a vertex $v\in V$ are defined as $d^-(v)=|\{(u,v)\in \vec{D}: u\in V\}|$ and $d^+(v)=|\{(v,u)\in \vec{D}: u\in V\}|$, respectively. The semidegree of $v$ is $d_0(v)=\min\{d^-(v), d^+(v)\}$. The minimum semidegree of the oriented graph $G=(V,\vec{D})$ is $\delta^0(G):=\min\{d^0(v):v\in V\}$. 

Analogously to undirected graphs, it is natural to ask for oriented graphs what minimum semidegree conditions would guarantee occurrence of long paths. The first case to consider are the directed paths, i.e. paths whose edges are directed the same way. Here Jackson~\cite{Ja} proved, as a corollary of a more general result, that $\delta^0(G)>k/2$ implies that $G$ contains a directed path with $k$ edges. The proof follows a path extension argument, similar to the textbook proof of Dirac's theorem~\cite{Di}. Subsequent years saw a large body of research on Hamilton cycles in oriented and directed graphs~\cite{Gh,Hg, HgTh1, HgTh2, KevKO, Kel, KelKO, ThE}. 

Motivated by the above works, as well as by results on degree conditions for occurrence of trees in undirected and directed graphs (see \cite{St}),  Stein~\cite{St} conjectured the following generalization of Jackson's theorem~\cite{Ja}. 
\begin{conj}[Stein~\cite{St}]\label{conj:Stein}
Every oriented graph with $\delta^0(G)>k/2$ contains every orientation of a path with $k$ edges.    
\end{conj}
\noindent
While it is easy to see that for $\delta^0(G)\geq k$ we can find greedily every orientation of a $k$-edge path, at present we do not know to replace $k$ with any smaller value (let alone $k/2$), for the statement of Conjecture~\ref{conj:Stein} to hold universally. 

After the directed paths, the next most natural case to consider are the \emph{alternating} paths, that is, paths in which every vertex sees either just incoming or just outgoing edges (this case also provides a construction for the bound of $k/2$ in Conjecture~\ref{conj:Stein}: a $k/2$-blowup of a directed cycle). Klimo\v{s}ov\'a and Stein~\cite{KSt} proved that $\delta^0(G)>3k/4$ implies existence of an alternating path with $k$ edges, and Chen, Hou and Zhou~\cite{ChHZ} improved their bound to $2k/3$ (up to additive constants). For the special case when $k$ is linear in the number of vertices of the oriented graph, Stein and Zárate-Guerén \cite{StZa} proved an approximate upper bound $(1/2+o(1))k$.

Our goal here is to improve the general upper bound further to $5k/8$. Similarly to~\cite{ChHZ} and~\cite{KSt}, we shall prove a slightly more general result, using the \emph{minimum pseudo-semidegree} $\bar{\delta}^0(G)$, which is defined, for an oriented graph without isolated vertices, as the minimum over all \emph{positive} vertex in- and out-degrees.

For better proof clarity we decided not to optimize additive constants. Consequently, we find it more convenient to deal with paths with $k$ vertices rather than edges.

\begin{theorem}\label{thm:main}
Every oriented graph $G$ with $\bar{\delta}^0(G)>5k/8$ contains an alternating path with $k$ vertices.
\end{theorem}

Stein and Z\'arate-Guer\'en~\cite{StZa}, using an elegant tensoring trick, observed that a digraph $H$ with more than $(\ell-1)|V(H)|$ edges contains a subgraph of minimum pseudo-semidegree at least $\ell/2$. Hence, as a corollary of Theorem~\ref{thm:main} we obtain the following improvement towards a conjecture of Addario-Berry, Havet, Linhares Sales, Reed and Thomass\'e~\cite{AHLRT}.

\begin{cor}
Every oriented graph $G$ with more than $(5k+4)|V(G)|/4$ edges contains an alternating path with $k$ vertices.
\end{cor}

\section{The proof of Theorem \ref{thm:main}.}
We start with some notation. Let $G=(V,\vec{D})$ be an oriented graph. For a vertex $v\in V$ and a vertex set $X\subseteq V$ let $N^+_{X}(v)$ and $N^-_{X}(v)$ denote $\{x\in X: (v,x)\in \vec{D}\}$ and $\{x\in X: (x,v)\in \vec{D}\}$ respectively. Let $d^+_{X}(v):=|N^+_{X}(v)|$ and $d^-_{X}(v):=|N^-_{X}(v)|$. Let $N^+(v):=N^+_{V}(v)$ and $N^-(v):=N^-_{V}(v)$ denote the out- and in-neighbourhoods of $v$.

Suppose now that $G$ is an oriented graph with $\bar{\delta}^0(G)>5k/8$. Let $P$ a maximum length alternating path in $G$, and assume for contradiction that $P$ has strictly fewer than $k$ vertices. We recall a relevant result of Klimo\v{s}ov\'a and Stein~\cite{KSt}, restated below in a convenient form. 

\begin{prop}[\cite{KSt}, Lemma 6]\label{prop:oddcase}
Let $k'\in \mathbb{N}$ and let $G$ be an oriented graph with $2\bar{\delta}^0(G)\geq k'$. Let $P$ be a longest alternating path in $G$ and suppose that $|P|$ is odd. Then $|P|\geq k' -1 $.    
\end{prop}
\noindent
Since we have $2\bar{\delta}^0(G) > 5k/4 \geq \lceil 5k/4 \rceil$, $|P|$ being odd would imply $|P|\geq \lceil 5k/4 \rceil-1\geq k$, a contradiction.

So, we may assume that $P$ has an even number of vertices $2m<k$. Let the vertices of $P$ be labeled $p_1,\dots ,p_{2m}$, in the order of the underlying undirected path. Without loss of generality assume that the edge $p_1p_2$ is oriented toward $p_2$. Let $O$ and $E$ be the sets of the odd and even indexed vertices $p_i$, respectively; we shall call them simply \emph{odd} and \emph{even} vertices. Note that every edge of $P$ is directed from an odd to an even vertex. Therefore, by definition of $\bar{\delta}^0(G)$, we have $d^+(v)>5k/8$ for every $v\in O$ and $d^-(w)>5k/8$ for every $w\in E$.

\begin{lemma}\label{lem:evenham}
There is an alternating spanning cycle on $V(P)=E\cup O$ such that each of its edges is directed from odd to even. 
\end{lemma}

\begin{proof}

We call an alternating path $P'\subseteq G$ \emph{respectable} if $V(P')=E\cup O$ and each edge of $P'$ is oriented from an odd to an even vertex (in particular, one end vertex of $P'$ is odd and one is even). The set of respectable paths is not empty, since $P$ is respectable. Call an even vertex $v\in E$ \emph{terminal} if is an end vertex of a respectable path, and let $T$ denote the set of terminal vertices. Call an odd vertex $v\in O$ \emph{starting} if it is an initial vertex of a respectable path, and let $S$ denote the set of starting vertices. So, for instance, $p_{2m}\in T$ and $p_1\in S$, as witnessed by a respectable path $P$.

Suppose $v\in S$ is a starting vertex of some respectable path $P'$. If there is an edge $(v,w)\in\vec{D}$ from some $w\notin E\cup O$, we would be able to extend $P'$ via the edge $(v,w)$ obtaining a longer alternating path. Thus we may assume that $N^{+}(v) \subseteq E\cup O$. 
So, since $d^+_{E\cup O}(v)=d^+(v)> 5k/8> (5/8)\cdot 2m$ and $|E|=|O|=m$, we must have $d_E^+(v)>0$ and also $d_O^+(v)>0$.
Similarly, for every terminal $w\in T$ we have $N^-(w)\subseteq E\cup O$, and in particular $d^-_O(w)>0$ and $d^-_E(w)>0$. 

Let now $a\in S$ be a starting vertex with $d^+_E(a)=\max\{d^+_E(v):v\in S\}$. 
Fix some respectable path $R$ starting at $a$, and denote its vertices $
a=r_1,r_2,\dots ,r_{2m}$ in the order of the underlying undirected path. Note that for each $r_i\in N^+_E(a)$ we have $r_{i-1}\in S$, since we have a respectable path $r_{i-1},\dots ,r_1,r_i,\dots,r_{2m}$. So, letting $A:=\{r_{i-1}\colon r_i \in N^+_E(a)\}$, we have $A\subseteq S$ and $|A|=d_E^+(a)$. 

\noindent
Now, by maximality of $d^+_E(a)$, for each vertex $u\in A$ we have 
$$d_E^+(u)\leq d_E^+(a)=|A|,$$ and therefore 
$$d_O^+(u)\geq \frac{5}{8}k-|A|> \frac{5}{8}\cdot 2m-|A|.$$ 
Let $s:=\sum_{u\in A}d_O^+(u)$. 
By the above, we have
$$s> \frac{5}{8} |A|\cdot 2m-|A|^2.$$ 
On the other hand, counting the edges via out-degrees (note that there are at most $\binom{|A|}{2}\leq \frac{|A|^2}{2}$ edges among the vertices in $A$, as $G$ is an oriented graph) and using $|O|=m$, yields  
$$s\leq \frac{|A|^2}{2}+|A|(m-|A|).$$
It follows that
$$\frac{|A|^2}{2}+|A|(m-|A|)> \frac{5}{8} |A|\cdot 2m-|A|^2,
$$
resulting in 
$$d^+_E(a)=|A|> \frac{m}{2}.$$

Now let $B$ be the set of all vertices $v\in T$ such that there exists a respectable path starting at $a$ and terminating at $v$; note that $B\neq \emptyset$, as $R$ is respectable and starts at $a$, so $r_{2m}\in B$. Let $b\in B$ be a vertex satisfying $d_O^-(b)=\max\{d_O^-(v):v\in B\}$.
Without loss of generality we may assume that $b=r_{2m}$. Note now that for each $r_i\in N_O^-(b)$ we have $r_{i+1}\in B$ due to the respectable path $a=r_1,\dots,r_i,r_{2m},\dots,r_{i+1}$. So, letting $C:=\{r_{i+1}\colon r_i\in N_O^-(b)\}$, we have $C\subseteq B$ and $|C|=d_O^-(b)$. 

\noindent
By maximality of $d_O^-(b)$, for each vertex $w\in C$ we have

$$d_O^-(w)\leq d_O^-(b)=|C|,
$$
and so
$$
d_E^-(w)\geq \frac{5}{8}k-|C|>\frac{5}{8}\cdot 2m-|C|.
$$
Let $t:=\sum_{w\in C} d_E^-(w)$. We obtain
$$t> |C|\cdot \frac{5}{8}\cdot 2m - |C|^2.
$$
On the other hand, counting the edges via in-degrees gives
$$t\leq \frac{|C|^2}{2}+|C|(m-|C|).
$$
Hence,
$$\frac{|C|^2}{2}+|C|(m-|C|)> |C|\cdot \frac{5}{8}\cdot 2m - |C|^2,
$$
which results in
$$d_O^-(b)=|C|> \frac{m}{2}.$$

Notice now that $|O|=m$, and we have shown that $N_O^-(b)\subseteq O$ and $\{r_i:r_{i+1}\in N_E^{+}(a)\}\subseteq O$ are of size greater than $m/2$. Therefore, by pigeonhole, for some $i$ we must have $r_i\in N_O^-(b)=N_O^-(r_{2m})$ and $r_{i+1}\in N^+_E(a)=N^+_E(r_1)$. This gives the desired cycle via the cyclic ordering $r_1,\dots r_i,r_{2m}, r_{2m-1},\dots,r_{i+1},r_1$. 
\end{proof}

We remark that Theorem~\ref{thm:main} can now be deduced directly from Lemma~\ref{lem:evenham} and Lemma 8 from~\cite{KSt}. Below we give an alternative argument that we believe carries some merit, providing a different perspective.

\begin{cor}\label{cor:evenodd}
We have $S=O$ and $T=E$. Consequently, for every $a\in O$ and $b\in E$ we have $N^+(a)\subseteq E\cup O$ and $N^-(b) \subseteq E\cup O$.
\end{cor} 

Let now $H$ be the simple bipartite graph obtained by taking all edges directed from odd to even vertices, and ignoring their directions. We claim that not only is $H$ hamiltonian (we already know this) but it satisfies, with room to spare, the Moon-Moser degree condition for hamiltonicity of bipartite graphs. To this end, let $d(v)$ denote the degree of a vertex $v\in E\cup O$ in the graph $H$. 

\begin{lemma}\label{lem:forgotten}
For every $1\leq \ell \leq m/2$ we have $$|\{v\in O: d(v)\leq \ell+1\}|<\ell.$$ And, symmetrically, $$|\{w\in E:d(w)\leq \ell+1\}<\ell.$$
\end{lemma}

\begin{proof}
Suppose for a contradiction that for some $1\leq \ell\leq m/2$ there exists a set $L\subseteq O$ (the statement for $E$ is proved in the same way) of $\ell$ vertices of $H$-degree at most $\ell+1$. Then each vertex $v\in L$ satisfies (using $d^+(v)>5k/8$ and $k\geq 2m+1$)
\begin{align*}   
d_O^+(v)&>\frac{5}{8} k-\ell-1\geq \frac{5}{8}(2m+1)+\frac{1}{8}-\ell-1\\
&=\frac{5}{4}m-\ell-\frac{1}{4}.
\end{align*}
It follows that 
$$|N^+(v)\cap L|> \frac{5}{4}m-\ell-\frac{1}{4} - (m-\ell)=\frac{m-1}{4}>\frac{\ell-1}{2},
$$
a contradiction, counting out-degrees inside $L$.
\end{proof}

\begin{lemma}\label{lem:MMcondition}
$G$ contains an alternating path $Q$ on $2m$ vertices, of the form $q_1\dots q_{2m}$, where (odd-indexed) $q_3,\dots,q_{2m-1}\in O$, (even-indexed) $q_4,\dots,q_{2m} \in E$, $q_2\in O$, $q_1\notin 
E\cup O$, and the edge $q_1q_2$ is directed from $q_1$ to $q_2$.
\end{lemma}

\begin{proof}
Recall that $d_O^+(v)>0$ for every $v \in O$. So, take two vertices $q_2,q_3\in O$ such that $(q_3,q_2)\in \vec{D}$.
Since $d^{-}(q_2)>0$, by the minimum pseudo-semidegree condition we have $d^{-}(q_2)>5k/8>(5/8)\cdot 2m$ in addition to $d^{+}(q_2)>(5/8)\cdot 2m$, as $q_2\in O$.
Since by Corollary~\ref{cor:evenodd} we have $N^+(q_2)\subseteq E\cup O$, and $G$ is an oriented graph, we conclude that $N^-(q_2)\setminus (E\cup O)\neq \emptyset$, and so $q_2$ has an in-neighbour $q_1 \notin E\cup O$. 

Let $O'=O\setminus \{q_2\}$, and let $E'$ be the set $E$ with an arbitrary vertex removed. Then, by Lemma~\ref{lem:forgotten}, the bipartite graph $H'=H[E'\cup O']$ satisfies the Moon-Moser condition (Corollary 1 in \cite{MM}): for any $1\leq \ell\leq |O'|/2=(m-1)/2$ the number of vertices in $O'$ of $H'$-degree at most $\ell$ (thus, $H$-degree at most $\ell+1$) is less than $\ell$, and symmetrically in $E'$. Therefore, $H'$ has a hamilton cycle, and, a fortiori, a hamilton path $q_3,\dots,q_{2m}$. This yields a desired alternating path $Q=q_1,q_2,q_3,\dots,q_{2m}$ in $G$.  
\end{proof}

Observe that Lemma~\ref{lem:evenham} applies to any alternating path on $2m$ vertices, and in particular, to $Q$, in place of $P$. Now consider the vertex $q_2$. Since $q_2\in O$, we have $N^+(q_2)\subset E\cup O$. On the other hand, $q_2$ plays the role of an `even' vertex in $Q$. 
So, by Lemma~\ref{lem:evenham} we must have $N^-(q_2)\subset E\cup O\cup \{q_1\}$. Put together, this would imply $d^-(q_2)+d^+(q_2)\leq 2m$, contradicting the fact that $d_0(q_2)> (5/8)\cdot 2m$, as $G$ is an oriented graph. This concludes the proof.

\end{document}